\documentclass[11pt]{llncs}

\usepackage{amssymb,amsmath,mathtools,url}
\usepackage[dvips]{graphics}

\usepackage{enumerate}
\usepackage{mleftright}
\usepackage[symbol]{footmisc}
\usepackage[affil-it]{authblk}
\usepackage{enumitem}   
\usepackage{url}

\usepackage{tabularx, booktabs}	
\newcolumntype{Y}{>{\centering\arraybackslash}X}

\makeatletter
\newcommand*{\myfnsymbolsingle}[1]{%
  \ensuremath{%
    \ifcase#1
    \or 
      *%
    \or 
      \dagger
    \or 
      \ddagger
    \or 
      \mathsection
    \or 
      \mathparagraph
    \else 
      \@ctrerr  
    \fi
  }%
}   
\makeatother

\usepackage{alphalph}
\newalphalph{\myfnsymbolmult}[mult]{\myfnsymbolsingle}{}

\newcommand{\F}{{\mathbb F}}

\newcommand{\ie}{i.e., }

\usepackage{vmargin}
\setpapersize{A4}
\setmarginsrb{1.75cm}{1.75cm}{1.75cm}{1.75cm}{0pt}{1.0cm}{0pt}{1.0cm}

\begin{document}

\title{Three Proofs of an Observation on Irreducible Polynomials over $\text{GF}(2)$}

\author{Robert Granger}
\institute{
Surrey Centre for Cyber Security\\
Department of Computer Science\\
University of Surrey\\
United Kingdom\\
\email{r.granger@surrey.ac.uk}
}

\maketitle

\begin{abstract}
We present three proofs of an observation of Ahmadi on the number of irreducible polynomials over $\text{GF}(2)$ with certain traces and cotraces, the most interesting of which uses an explicit natural bijection. We also present two proofs of a related observation.
\end{abstract}


\section{Introduction}\label{sec:intro}


For $n \ge 2$ let $\mathcal{I}_n$ denote the set of all irreducible degree $n$ polynomials in $\F_{2}[x]$, and let $\text{Tr}_n: \F_{2^n} \rightarrow \F_2: \alpha \mapsto \alpha + \alpha^2 + \alpha^{2^2} + \cdots + \alpha^{2^{n-1}}$ denote the absolute trace function. 
For a polynomial $f = x^n + f_{n-1}x^{n-1} + \cdots + f_1x + 1 \in \mathcal{I}_n$, if $\alpha$ is a root of $f$ then $f_{n-1} = \text{Tr}_n(\alpha)$ and $f_{1} = \text{Tr}_n(\alpha^{-1})$: $f_{n-1}$ and $f_{1}$ are known as the trace and cotrace respectively.
We partition $\mathcal{I}_n$ into four sets $S_{i,j}(n)$ with $i,j \in \F_2$ by placing each $f \in \mathcal{I}_n$ into $S_{f_{n-1},f_1}(n)$. 
Table~\ref{table1} contains the cardinality of these sets for $2 \le n \le 32$ (note that we do not define $S_{i,j}(1)$).
Elements of $S_{1,1}(n)$ are useful for practical applications since they give rise to representations of $\F_{2^{n \cdot 2^{l}}}$ for all $l \ge 1$ via the iteration of the so-called $Q$-transform~\cite{meyn} (cf. \S\ref{sec:parity}), provided that $n \ne 3$~\cite{niederreiter}.

It is clear that for any $n \ge 3$ the sets $S_{0,1}(n)$ and $S_{1,0}(n)$ have the same cardinality, 
since any member $f$ of one set can be mapped to a corresponding member of the other set via the reciprocal transform $f^*(x) = x^n f(1/x)$, which reverses the coefficients of $f$: since this transform is invertible (indeed it is its own inverse) it gives a natural bijection between the two sets, in the sense that it is simple and has explanatory power.
Ahmadi observed that for odd $n \ge 3$ the sets $S_{0,0}(n)$ and $S_{1,1}(n)$ also have the same cardinality~\cite{omran}, which raises the question of whether or not there exists a natural bijection between them, just as for $S_{0,1}(n)$ and $S_{1,0}(n)$? 

There exist bijective proofs of numerous combinatorial identities: indeed, Stanley has exhibited  $250$ `Bijective Proof Problems' of various levels of difficulty, including 27 open problems~\cite{Stanley}. Occasionally, a natural bijection can illuminate the relation between two sets of the same cardinality. One such example is Benjamin and Bennett's elegant solution~\cite{Dilcue} to a question posed by Corteel, Savage, Wilf and Zeilberger, which asked for a bijective explanation of the fact that among ordered pairs of polynomials of degree $n$ over $\F_{2}[x]$, there are as many coprime pairs as there are non-coprime pairs~\cite{CSWZ}. Benjamin and Bennett constructed such a bijection by applying Euclid's algorithm to any pair, flipping the final remainder bit and then reversing Euclid's algorithm using the same quotients. The main purpose of the present work is to exhibit a natural bijection which explains Ahmadi's observation.

The author further observed that for even $n$, the difference $|S_{1,1}(n)| - |S_{0,0}(n)|$ is equal to the number of trace $1$ irreducibles of degree $n/2$. 
Before presenting our bijective proof of Ahmadi's observation in \S\ref{sec:bij}, for good measure we first prove his and our observations in two different ways, in \S\ref{sec:easy} and~\S\ref{sec:nied}, each proof having its own merits. We finish by presenting a proposition on the parity of $|S_{1,1}(n)|$ in \S\ref{sec:parity}, which arises from similar considerations. For reference and clarity we now state our two main theorems explicitly.

\begin{theorem}\label{mainthm}
For odd $n \ge 3$ the sets $S_{0,0}(n)$ and $S_{1,1}(n)$ have the same cardinality.
\end{theorem}

\begin{table}[h!]
\caption{Cardinality of $S_{i,j}(n)$ for $2 \le n \le 32$}
\begin{center}\label{table1}
\begin{tabular}{c|c|c|c|c}
\hline
$n$ & $|S_{0,0}(n)|$ & $|S_{0,1}(n)|$ & $|S_{1,0}(n)|$ & $|S_{1,1}(n)|$\\
\hline
2 & 0 & 0 & 0 & 1 \\
3 & 0 & 1 & 1 & 0 \\
4 & 0 & 1 & 1 & 1 \\
5 & 2 & 1 & 1 & 2 \\
6 & 1 & 3 & 3 & 2 \\
7 & 4 & 5 & 5 & 4 \\
8 & 7 & 7 & 7 & 9 \\
9 & 14 & 14 & 14 & 14\\
10 & 21 & 27 & 27 & 24 \\
11 & 48 & 45 & 45 & 48 \\
12 & 81 & 84 & 84 & 86  \\
13 & 154 & 161 & 161 & 154 \\
14 & 285 & 291 & 291 & 294 \\
15 & 550 & 541 & 541 & 550 \\
16 & 1001 & 1031 & 1031 & 1017 \\
17 & 1926 & 1929 & 1929 & 1926 \\
18 & 3626 & 3626 & 3626 & 3654 \\
19 & 6888 & 6909 & 6909 & 6888 \\
20 & 13041 & 13122 & 13122 & 13092 \\
21 & 24998 & 24931 & 24931 & 24998 \\
22 & 47565 & 47667 & 47667 & 47658 \\
23 & 91124 & 91237 & 91237 & 91124 \\
24 & 174652 & 174698 & 174698 & 174822 \\
25 & 335588 & 335500 & 335500 & 335588 \\
26 & 644805 & 645435 & 645435 & 645120 \\
27 & 1242822 & 1242682 & 1242682 & 1242822 \\
28 & 2396385 & 2396520 & 2396520 & 2396970 \\
29 & 4627850 & 4628545 & 4628545 & 4627850 \\
30 & 8946665 & 8947923 & 8947923 & 8947756 \\
31 & 17319148 & 17317685 & 17317685 & 17319148 \\
32 & 33551833 & 33554983 & 33554983 & 33553881 \\
\hline
\end{tabular}
\end{center}
\end{table}

\begin{theorem}\label{evenn}
For even $n$, the difference $|S_{1,1}(n)| - |S_{0,0}(n)|$ is equal to the number of trace $1$ irreducibles of degree $n/2$. In particular, we have
\[
|S_{1,1}(n)| - |S_{0,0}(n)| = \frac{1}{n}\sum_{\substack{{d \mid n/2} \\
                              d \ \text{odd}}} \mu(d) 2^{n/2d},  
\]
where $\mu(\cdot)$ is the M\"obius function, which is defined by:
\[
\mu(n) = \begin{cases}
1 & \text{if } n = 1,\\
(-1)^k & \text{if } n \text{ is a product of $k$ distinct primes},\\ 
0 & \text{otherwise}.
\end{cases}
\]
\end{theorem}

\section{First proofs of the observations}\label{sec:easy}

Our first proofs of Theorems~\ref{mainthm} and~\ref{evenn} are easy and the most direct, but are perhaps the least illuminating since they use two well-known theorems.
\vspace{3mm}
\newline
\noindent {\em First proof of Theorem~\ref{mainthm}.}
For $n \ge 1$ it is well known that 
\[
|\mathcal{I}_n| = \frac{1}{n} \sum_{d \mid n} \mu(d) 2^{n/d}.
\]
It is also well known (see~\cite{carlitz}) that the number of binary irreducibles of degree $n \ge 1$ with trace $1$ is 
\begin{equation}\label{trace1count}
\frac{1}{2n} \sum_{\substack{{d \mid n} \\
                              d \ \text{odd}}} \mu(d) 2^{n/d}.
\end{equation}
Assume now that $n \ge 2$. Then~(\ref{trace1count}) equals $|S_{1,1}(n)| + |S_{1,0}(n)|$. Furthermore, since 
$|S_{1,0}(n)| = |S_{0,1}(n)|$ we have
\begin{eqnarray}
\nonumber |S_{1,1}(n)| - |S_{0,0}(n)| &=& (|S_{1,1}(n)| + |S_{1,0}(n)|) - (|S_{0,1}(n)| + |S_{0,0}(n)|)\\
\nonumber                       &=& (|S_{1,1}(n)| + |S_{1,0}(n)|) - (|\mathcal{I}_n| - |S_{1,1}(n)| - |S_{1,0}(n)|)\\
\nonumber                       &=& 2(|S_{1,1}(n)| + |S_{1,0}(n)|) - |\mathcal{I}_n|\\
\label{terms}                  &=& \frac{1}{n} \sum_{\substack{{d \mid n} \\
                              d \ \text{odd}}} \mu(d) 2^{n/d} - \frac{1}{n} \sum_{d \mid n} \mu(d) 2^{n/d} = -\frac{1}{n}\sum_{\substack{{d \mid n} \\
                              d \ \text{even}}} \mu(d) 2^{n/d}  
\end{eqnarray}
If $n$ is odd then the final sum in expression~(\ref{terms}) is empty, which proves Theorem~\ref{mainthm}. \qed
\vspace{3mm}

\noindent {\em First proof of Theorem~\ref{evenn}.}
From~(\ref{terms}) we have
\begin{eqnarray}
\label{line3} |S_{1,1}(n)| - |S_{0,0}(n)| &=& -\frac{1}{n}\sum_{\substack{{d \mid n} \\
                              d \ \text{even}}} \mu(d) 2^{n/d} \\
\label{line4} &=& -\frac{1}{n}\sum_{\substack{{2d \mid n} \\
                              d \ \text{odd}}} \mu(2d) 2^{n/2d} \\
\label{line5} &=& \frac{1}{n}\sum_{\substack{{d \mid n/2} \\
                              d \ \text{odd}}} \mu(d) 2^{n/2d},               
\end{eqnarray}                  
where expression~(\ref{line4}) discounts all those $d$ in the sum in equation~(\ref{line3}) which are divisible by $4$, since
$\mu$ of such $d$ is zero. The sum in~(\ref{line5}) is nothing but~(\ref{trace1count}) but for argument $n/2$, as claimed. \qed

\section{Second proofs of the observations}\label{sec:nied}

Our second proofs of Theorems~\ref{mainthm} and~\ref{evenn} are based on Niederreiter's explicit count of $|S_{1,1}(n)|$~\cite{niederreiter} and arguably give more insight than our first proofs.
\vspace{3mm}
\newline
\noindent {\em Second proof of Theorem~\ref{mainthm}.}
Let $n \ge 1$. For $i \in \F_2$ let $N_i(n) = \#\{ \alpha \in \F_{2^n}^{\times} \mid \text Tr_n(\alpha) = \text{Tr}_n(\alpha^{-1}) = i\}$. Niederreiter expressed $N_1(n)$ as follows:
\begin{eqnarray}
\nonumber N_1(n) &=& \sum_{\alpha \in \F_{2^n}^{\times}} \bigg( \frac{1}{2} \sum_{a \in \F_2} (-1)^{a(\text{Tr}_n(\alpha) + 1)} \bigg) 
                                                                           \bigg( \frac{1}{2} \sum_{b \in \F_2} (-1)^{b(\text{Tr}_n(\alpha) + 1)} \bigg)\\
\nonumber            &=& \frac{1}{4} \sum_{a,b \in \F_2} (-1)^{a+b} \sum_{\alpha \in \F_{2^n}^{\times}} (-1)^{\text{Tr}_{n}(a \alpha + b \alpha^{-1})} \\
\label{kloostersum}            &=& \frac{1}{4} \big((2^n - 1) +1 + 1 + \sum_{\alpha \in \F_{2^n}^{\times}} (-1)^{\text{Tr}_{n}(\alpha + \alpha^{-1})} \big),
\end{eqnarray}
where in~(\ref{kloostersum}) the $2^n - 1$ corresponds to $a = b = 0$, the two $+1$'s correspond to $a \ne b$, and the final sum
corresponds to $a = b = 1$. Note that the final sum is the well-known Kloosterman sum evaluated at $1$, but for our purposes we need not evaluate it. In particular, using a similar argument we have:
\begin{eqnarray}
\nonumber N_0(n) &=& \sum_{\alpha \in \F_{2^n}^{\times}} \bigg( \frac{1}{2} \sum_{a \in \F_2} (-1)^{a \text{Tr}_n(\alpha)} \bigg) 
                                                                           \bigg( \frac{1}{2} \sum_{b \in \F_2} (-1)^{b \text{Tr}_n(\alpha)} \bigg)\\
\nonumber            &=& \frac{1}{4} \sum_{a,b \in \F_2}  \sum_{\alpha \in \F_{2^n}^{\times}} (-1)^{\text{Tr}_{n}(a \alpha + b \alpha^{-1})} \\
\nonumber           &=& \frac{1}{4} \big((2^n - 1) -1 - 1 + \sum_{\alpha \in \F_{2^n}^{\times}} (-1)^{\text{Tr}_{n}(\alpha + \alpha^{-1})} \big),
\end{eqnarray}
and therefore $N_0(n) = N_1(n) - 1$.

Now let $G_i(n) = \#\{ \alpha \in \F_{2^n}^{\times} \mid \text Tr_n(\alpha) = \text{Tr}_n(\alpha^{-1}) = i \ \text{and} \ \F_2(\alpha) = \F_{2^n} \}$,
\ie $G_i(n)$ is the cardinality of the subset of elements counted by $N_i(n)$ which are roots of irreducible degree $n$ polynomials. Since any irreducible degree $n$ polynomial has precisely $n$ roots in $\F_{2^n}$ we see that $|S_{i,i}(n)| = \frac{1}{n} G_i(n)$ for $n \ge 2$. For each $\alpha \in \F_{2^n}^{\times}$ there is a uniquely determined irreducible polynomial in $\F_2[x]$ of degree $d \mid n$ for which $\alpha$ is a root, and so by transitivity of the trace, for this $d$ we have
\begin{equation}\label{transitivity}
\text{Tr}_n(\alpha) = \text{Tr}_d \big( \frac{n}{d} \alpha \big) = \frac{n}{d} \text{Tr}_d(\alpha).
\end{equation}
Thus $\text{Tr}_n(\alpha) = 1$ if and only if $n/d$ is odd and $\text{Tr}_d(\alpha) = 1$, and likewise for $\alpha^{-1}$. 
Niederreiter therefore deduces that
\begin{equation}\label{N1}
N_1(n) = \sum_{\substack{{d \mid n} \\
                              n/d \ \text{odd}}} G_1(d).
\end{equation}
Also, by~(\ref{transitivity}) $\text{Tr}_n(\alpha) = 0$ if and only if $n/d$ is even and it does not matter what $\text{Tr}_d(\alpha)$ is, or 
$n/d$ is odd and $\text{Tr}_d(\alpha) = 0$, and likewise for $\alpha^{-1}$.
In the former case the contribution to $N_0(n)$ is simply the cardinality of the largest subfield of $\F_{2^n}$ such that $n/d$ is
even, minus $1$ since the zero element is not counted by $N_0(n)$. We therefore deduce that
\begin{equation}\label{N0}
N_0(n) = |\F_{2^{\text{max} \{ d \ \mid \  n/d \ \text{even} \} }}^{\times} | + \sum_{\substack{{d \mid n} \\
                              n/d \ \text{odd}}} G_0(d) .
\end{equation}
If $n$ is odd then by M\"obius inversion~(\ref{N1}) gives
\[
G_1(n) = \sum_{d \mid n} \mu(d) N_1(n/d),
\]
while the first term in~(\ref{N0}) becomes empty and M\"obius inversion gives
\begin{eqnarray}
\nonumber G_0(n) &=& \sum_{d \mid n} \mu(d) N_0(n/d) \\
\nonumber            &=& \sum_{d \mid n} \mu(d) (N_1(n/d) - 1) \\
\nonumber            &=& \sum_{d \mid n} \mu(d) N_1(n/d)  - \sum_{d \mid n} \mu(d) \\
\nonumber            &=& G_1(n) \ \text{for} \ n \ge 2.
\end{eqnarray}
Since $|S_{i,i}(n)| = \frac{1}{n} G_i(n)$ for $n \ge 2$ we have $|S_{0,0}(n)| = |S_{1,1}(n)|$ for odd $n \ge 3$ which reproves Theorem~\ref{mainthm}. \qed
\vspace{3mm}

\noindent {\em Second proof of Theorem~\ref{evenn}.}
Let $k \ge 1$ and let $t \ge1$ be odd. Then for argument $2^k t$, 
the $d$ occurring in the sums of~(\ref{N1}) and~(\ref{N0}) are of the form $d = 2^k e$ with $e$ a positive divisor of $t$. We thus have
\begin{equation}\label{N1even}
N_1(2^k t) = \sum_{e \mid t} G_1(2^k e),
\end{equation}
and
\begin{equation}\label{N0even}
N_0(2^k t) = | \F_{2^{2^{k-1} t}}^{\times} | + \sum_{e \mid t} G_0(2^k e).
\end{equation}
Since one cannot immediately apply M\"obius inversion to~(\ref{N1even}) and~(\ref{N0even}) to obtain $G_1(2^k t)$ and 
$G_0(2^k t)$, for any integer $m \ge 1$ and $i = 0,1$ define 
\[
H_i(m) = \sum_{d|m} G_i(2^k d).
\]
Then by M\"obius inversion we have
\begin{equation}\label{Gexp}
G_i(2^k m) = \sum_{d \mid m} \mu \bigg(\frac{m}{d} \bigg) H_i(d).
\end{equation}
If $m$ is odd then by the definition of $H_i$ and by~(\ref{N1even}) and~(\ref{N0even}) respectively, we have 
$H_1(m) = N_1(2^k m)$ and $H_0(m) = N_0(2^k m) - | \F_{2^{2^{k-1} m}}^{\times} |$. Let $n = 2^k m$ with $k \ge 1$ and $m \ge 1$ odd. Then rewriting~(\ref{Gexp}) using these equations respectively, we obtain
\[
G_1(2^k m) = \sum_{d \mid m} \mu \bigg( \frac{m}{d} \bigg) N_1(2^k d),
\]
and
\begin{eqnarray}
\nonumber G_0(2^k m) &=& \sum_{d \mid m} \mu \bigg( \frac{m}{d} \bigg)  (N_0(2^k d) - | \F_{2^{2^{k-1} d}}^{\times} |) \\
\nonumber               &=& \sum_{d \mid m} \mu \bigg( \frac{m}{d} \bigg)  N_0(2^k d)  - \sum_{d \mid m} \mu \bigg( \frac{m}{d} \bigg) (2^{nd/2m} - 1) \\
\nonumber       &=& \sum_{d \mid m} \mu \bigg( \frac{m}{d} \bigg)  (N_1(2^k d) - 1) - \sum_{d \mid m} \mu ( d) (2^{n/2d}-1) \\
\label{final2}                  &=& G_1(2^k m) - \sum_{d \mid m} \mu ( d) 2^{n/2d}.
\end{eqnarray}
Dividing equation~(\ref{final2}) by $n$ reproves Theorem~\ref{evenn}. \qed

\section{An explicit bijection between $S_{0,0}(n)$ and $S_{1,1}(n)$ for odd $n$}\label{sec:bij}

We now present a bijective proof of Theorem~\ref{mainthm}. Crucial to our bijection are the following two transforms. Let
$\psi : \mathcal{I}_n \rightarrow \mathcal{I}_n : f \mapsto (x+1)^n f(\frac{1}{x+1})$, which has inverse
$\psi^{-1} : f \mapsto x^n f(\frac{x+1}{x})$, as is easily verified.
Since the arguments of $f$ in $\psi$ and $\psi^{-1}$ are invertible fractional linear transformations, they map irreducibles to irreducibles and are thus well-defined. 

We observed that under $\psi$ and $\psi^{-1}$, which $S_{i',j'}(n)$ an element of $S_{i,j}(n)$ maps to depends only on $i$ and $j$, the parity of $n$ and the parity of the number of monomials $x^k$ in the range $2 \le k \le n-2$ which have odd exponent, the latter of which motivates the following equivalent definition.

\begin{definition}
For a polynomial $f \in \mathcal{I}_n$ we define its signature $\sigma_f \in \F_2$ to be $\sum_{k = 2}^{n-2} k f_k \pmod{2}$.
\end{definition}

We have the following important lemma.

\begin{lemma}\label{lemma2.1}
Let $n \ge 3$ be odd, let $f \in S_{1,1}(n)$ and let $g \in S_{0,0}(n)$. Then 
\begin{enumerate}[label={(\roman*)}]
\item \[
\psi(f) \in \begin{cases}
S_{0,0}(n) &\mbox{if } \ \sigma_f  = 1 \\ 
S_{0,1}(n) & \mbox{if } \ \sigma_f = 0 \end{cases}
\]
\item \[
\psi^{-1}(f) \in \begin{cases}
S_{0,0}(n) &\mbox{if } \ \sigma_f  = 0 \\ 
S_{1,0}(n) & \mbox{if } \ \sigma_f = 1 \end{cases}
\]
\item \[
\psi(g) \in \begin{cases}
S_{1,1}(n) &\mbox{if } \ \sigma_g  = 1 \\ 
S_{1,0}(n) & \mbox{if } \ \sigma_g = 0 \end{cases}
\]
\item \[
\psi^{-1}(g) \in \begin{cases}
S_{1,1}(n) &\mbox{if } \ \sigma_g  = 0 \\ 
S_{0,1}(n) & \mbox{if } \ \sigma_g = 1 \end{cases}
\]
\end{enumerate}
\end{lemma}

\begin{proof}
For part (i), observe that
\begin{equation}\label{eqn1}
\psi(f) = (x+1)^n + (x+1)^{n-1} + \sum_{k=2}^{n-2} f_k (x + 1)^{n-k} + (x+1) + 1
\end{equation}
The coefficient of $x^{n-1}$ in~(\ref{eqn1}) is ${n \choose n-1} + 1 = n + 1 \equiv 0 \pmod{2}$, since $n$ is odd. The coefficient of $x$ in~(\ref{eqn1}) is 
\[
{n \choose 1} + {n-1 \choose 1} + \sum_{k = 2}^{n-2} f_k {n-k \choose 1} + 1 = n + (n - 1) + 
\sum_{k = 2}^{n-2} f_k (n-k) +1 \equiv 
\sum_{k=2}^{n-2} f_k + \sum_{k = 2}^{n-2} k f_k \pmod{2}.
\]
Since all irreducibles polynomials in $\F_2[x]$ of degree $ > 1$ necessarily have an odd number of terms (for otherwise $x+1$ would be a factor), we have $\sum_{k=2}^{n-2} f_k \equiv 1 \pmod{2}$, which completes the proof of part (i). For part (ii) 
observe that 
\begin{equation}\label{eqn2}
\psi^{-1}(f) = (x+1)^n + (x+1)^{n-1}x + \sum_{k=2}^{n-2} f_k (x + 1)^k x^{n-k} + (x+1)x^{n-1} + x^n.
\end{equation}
The coefficient of $x^{n-1}$ in~(\ref{eqn2}) is
\[
{n \choose n-1} + {n-1 \choose n-2} + \sum_{k = 2}^{n-2} f_k {k \choose k-1} + 1 = n + (n - 1) + \sum_{k = 2}^{n-2} f_k k + 1 \equiv 
\sum_{k = 2}^{n-2} k f_k \pmod{2}.
\]
The coefficient of $x$ in~(\ref{eqn2}) is ${n \choose 1} + 1 = n + 1 \equiv 0 \pmod{2}$. This completes the proof of
part (ii). For part (iii) observe that
\begin{equation}\label{eqn3}
\psi(g) = (x+1)^n + \sum_{k=2}^{n-2} g_k (x + 1)^{n-k} + 1.
\end{equation}
The coefficient of $x^{n-1}$ in~(\ref{eqn3}) is ${n \choose n-1} = n \equiv 1 \pmod{2}$, since $n$ is odd. The coefficient of $x$ in~(\ref{eqn3}) is 
\[
{n \choose 1} + \sum_{k = 2}^{n-2} g_k {n-k \choose 1} = n + \sum_{k = 2}^{n-2} g_k (n-k) \equiv 
1 + \sum_{k=2}^{n-2} g_k + \sum_{k = 2}^{n-2} k g_k \equiv \sum_{k = 2}^{n-2} k g_k \pmod{2},
\]
which proves part (iii). For part (iv) observe that
\begin{equation}\label{eqn4}
\psi^{-1}(g) = (x+1)^n + \sum_{k=2}^{n-2} g_k (x + 1)^k x^{n-k} + x^n.
\end{equation}
The coefficient of $x^{n-1}$ in~(\ref{eqn4}) is
\[
{n \choose n-1} + \sum_{k = 2}^{n-2} g_k {k \choose k-1} = n + \sum_{k = 2}^{n-2} g_k k \equiv 
1 + \sum_{k = 2}^{n-2} k g_k \pmod{2}.
\]
The coefficient of $x$ in~(\ref{eqn4}) is ${n \choose 1} = n \equiv 1 \pmod{2}$, which completes the proof of part (iv)
and the lemma. \qed
\end{proof}

We now reprove Theorem~\ref{mainthm} with an explicit bijection.
\vspace{3mm}
\newline
\noindent {\em Third proof of Theorem~\ref{mainthm}.}
Let $f \in S_{1,1}(n)$ and define a map 
$\phi : S_{1,1}(n) \rightarrow S_{0,0}(n)$ by 
\[
\phi(f) := \begin{cases}
\psi(f) &\mbox{if } \ \sigma_f = 1 \\ 
\psi^{-1}(f) & \mbox{if } \ \sigma_f = 0. \end{cases}
\]
Also, let $g \in S_{0,0}(n)$ and define a map 
$\rho: S_{0,0}(n) \rightarrow S_{1,1}(n)$ by 
\[
\rho(g) := \begin{cases}
\psi(g) &\mbox{if } \ \sigma_g = 1 \\ 
\psi^{-1}(g) & \mbox{if } \ \sigma_g = 0. \end{cases}
\]
We will show that $\phi$ and $\rho$ are inverse to one another. 
Firstly, if $\sigma_f = 1$ then by Lemma~\ref{lemma2.1}(i) we have $\phi(f) = \psi(f) \in S_{0,0}(n)$.
Since $\psi^{-1}(\psi(f)) = f \in S_{1,1}(n)$, by Lemma~\ref{lemma2.1}(iv) we must have $\sigma_{\psi(f)} = 0$.
Hence $\rho(\phi(f)) = f$ in this case. Furthermore, if $\sigma_f = 0$ then by Lemma~\ref{lemma2.1}(ii) we have 
$\phi(f) = \psi^{-1}(f) \in S_{0,0}(n)$. Since $\psi(\psi^{-1}(f)) = f \in S_{1,1}(n)$, by Lemma~\ref{lemma2.1}(iii) we must have $\sigma_{\psi^{-1}(f)} = 1$. Hence $\rho(\phi(f)) = f$ in this case too and $\rho$ is a left inverse for $\phi$.

Secondly, if $\sigma_g$ = 1 then by Lemma~\ref{lemma2.1}(iii) we have $\rho(g) = \psi(g) \in S_{1,1}(n)$.
Since $\psi^{-1}(\psi(g)) = g \in S_{0,0}(n)$, by Lemma~\ref{lemma2.1}(ii) we must have $\sigma_{\psi(g)} = 0$.
Hence $\phi(\rho(g)) = g$ in this case. Furthermore, if $\sigma_g = 0$ then by Lemma~\ref{lemma2.1}(iv) we have 
$\rho(g) = \psi^{-1}(g) \in S_{1,1}(n)$. Since $\psi(\psi^{-1}(g)) = g \in S_{0,0}(n)$, by Lemma~\ref{lemma2.1}(i) we must have 
$\sigma_{\psi^{-1}(g)} = 1$. Hence $\phi(\rho(g)) = g$ in this case too and $\rho$ is a right inverse for $\phi$. Thus $\phi$ and $\rho$ are inverse to one another. \qed

\vspace{3mm}
\subsection{An open problem for even $n$}

To complement the above proof of Theorem~\ref{mainthm}, it would be desirable to have a bijective proof of 
Theorem~\ref{evenn}, \ie a natural map between $S_{1,1}(n)$ and $S_{0,0}(n)$ union the set of trace $1$ irreducibles of 
degree $n/2$, when $n$ is even. One obstruction however is that the subset of $S_{0,0}(n)$ consisting of elements with signature $0$ maps to itself under the action on $\mathcal{I}_n$ of the group generated by the reciprocal transform and $\psi$, which is isomorphic to $GL_2(\F_2)$ (see~\cite{michon} for a classification of this action). Similarly, the subset of 
$S_{1,1}(n)$ consisting of elements with signature $1$ maps to itself under this action. Hence, if there exists such a bijection then other more sophisticated maps will be required. One possible approach consists of first factoring members of $S_{0,0}(n)$ and $S_{1,1}(n)$ over $\F_4$ into two degree $n/2$ 
(conjugate) irreducibles and acting on either factor by carefully chosen elements of $GL_2(\F_4)$ according to some arithmetic characteristics, just as we did with $\psi$ and $\psi^{-1}$, since then all polynomials concerned are of the same degree. However, the details of this action are naturally more complicated than the one arising from $GL_2(\F_2)$ and we leave its study and finding an explicit bijection as an open problem.

\section{The parity of $|S_{1,1}(n)|$}\label{sec:parity}

In this final short section we present an elementary result whose proof arises from simple transforms of polynomials and bijections. 
We first recall some relevant definitions and supporting results.

A polynomial $f \in \F_2[x]$ is said to be {\em self-reciprocal} if $f^* = f$. Let the set of degree $n$ self-reciprocal irreducible (SRI) polynomials in $\F_2[x]$ with trace $1$ be denoted by $\text{SRI}_{1}(n)$. For a degree $n$ polynomial $f$ the $Q$-transform of $f$, denoted $f^Q$, is defined to be $x^n f(x + 1/x)$, which is self-reciprocal and of degree $2n$.
A useful and well-known result -- originally due to Varshamov and Garakov~\cite{varshamov} and later generalised by 
Meyn~\cite{meyn} -- is that $f^Q$ is irreducible if and only if $f$ is irreducible and $f_1 = 1$. We have the following proposition.

\begin{proposition}
$|S_{1,1}(n)| \equiv 1 \pmod{2}$ if and only if $n = 2^k$ with $k \ge 1$.
\end{proposition}

\begin{proof}
The reciprocal transform acts on $S_{1,1}(n)$, partitioning it into pairs of distinct polynomials $(f,f^*)$ and a set of fixed points, namely 
$\text{SRI}_{1}(n)$. Hence $|S_{1,1}(n)| \equiv |\text{SRI}_{1}(n)| \pmod{2}$ and we need only determine the parity of 
$|\text{SRI}_{1}(n)|$. For odd $n > 1$ there are no SRIs, since if $\alpha$ is a root of an SRI then so is $1/\alpha$, and so 
the number of roots must be even. Therefore let $n$ be even. For any $f \in \text{SRI}_{1}(n)$ there 
exists a unique $f'$ of degree $n/2$ such that $f = f'^Q$ (see for instance~\cite[Lemma 6]{ahmadivega}). 
One may thus partition $\text{SRI}_{1}(n)$ into pairs of distinct polynomials $(f, (f'^*)^Q)$ and a set of fixed points for 
which $f' = f'^{*}$.
These fixed points are precisely $\text{SRI}_{1}(n/2)$, since by the Varshamov-Garakov criterion $f'$ is irreducible and 
$f_{1}^{'} = 1$, and the trace $f_{n/2 - 1}^{'}$ equals $f_{1}^{'}$.
Hence $|\text{SRI}_{1}(n)| \equiv  |\text{SRI}_{1}(n/2)|  \pmod{2}$. If $n = 2^k m$ with odd $m > 1$, then applying this
descent step repeatedly gives
\[
|S_{1,1}(2^k m)| \equiv |\text{SRI}_{1}(2^k m)| \equiv |\text{SRI}_{1}(2^{k-1} m)| \equiv \cdots \equiv |\text{SRI}_{1}(m)| 
\equiv 0 \pmod{2}.
\] 
On the other hand, if $n = 2^k$ then descending as before gives $|S_{1,1}(2^k)| \equiv |\text{SRI}_{1}(1)| \pmod{2}$. Since $x + 1$ is the only element of $\text{SRI}_1(1)$, the result follows. \qed
\end{proof}

Note that one could in principle analyse Niederreiter's (complicated) explicit formulae~\cite{niederreiter} for $|S_{1,1}(n)|$ in order to obtain this result. However, the above approach is perhaps more enlightening. 

\section*{Acknowledgements}
This work was supported by the Engineering and Physical Sciences Research Council via grant number EP/W021633/1. I would like to thank Omran Ahmadi for informing me of his observation, which motivated the search for a bijective proof of Theorem~\ref{mainthm}.



\bibliographystyle{plain}
\bibliography{3proofsFINALbib}

\end{document}